\documentclass[sn-mathphys-num]{sn-jnl}% Math and Physical Sciences Numbered Reference Style 
\usepackage{graphicx}%
\usepackage{multirow}%
\usepackage{amsmath,amssymb,amsfonts}%
\usepackage{amsthm}%
\usepackage{mathrsfs}%
\usepackage[title]{appendix}%
\usepackage{xcolor}%
\usepackage{textcomp}%
\usepackage{manyfoot}%
\usepackage{booktabs}%
\usepackage{algorithm}%
\usepackage{algorithmicx}%
\usepackage{algpseudocode}%
\usepackage{listings}%
\usepackage{adjustbox}
\usepackage{tabularx}
\theoremstyle{thmstyleone}%
\theoremstyle{thmstyletwo}%

\theoremstyle{thmstylethree}%

\usepackage{url}
\usepackage{footmisc}
\begin{document}

\title{Spatial Optimization of Autonomous Vehicle Assignment Based on Distance-Driven Demand and Customer Patience}

%%=============================================================%%
%% GivenName	-> \fnm{Joergen W.}
%% Particle	-> \spfx{van der} -> surname prefix
%% FamilyName	-> \sur{Ploeg}
%% Suffix	-> \sfx{IV}
%% \author*[1,2]{\fnm{Joergen W.} \spfx{van der} \sur{Ploeg} 
%%  \sfx{IV}}\email{iauthor@gmail.com}
%%=============================================================%%

\author[1]{\fnm{Niloufar} \sur{Mirzavand Boroujeni}}\email{mirza047@umn.edu}

\author[1]{\fnm{Nasim} \sur{Mirzavand Boroujeni}}\email{mirza048@umn.edu}

%\equalcont{These authors contributed equally to this work.}

\author*[2]{\fnm{Nima} \sur{Moradi}}\email{nima.moradi@mail.concordia.ca}

\author[3]{\fnm{Saeed} \sur{Jamalzadeh}}\email{s.jamalzadeh@ou.edu}

%\equalcont{These authors contributed equally to this work.}

\affil[1]{\orgdiv{Department of Industrial and Systems Engineering}, \orgname{University of Minnesota}, \orgaddress{%\street{Street}, 
 \city{Minneapolis}, \postcode{55455}, \state{MN}, \country{United States}}}

\affil*[2]{\orgdiv{Information Systems Engineering}, \orgname{Concordia University}, \orgaddress{\street{1455 Blvd. De Maisonneuve Ouest}, \city{Montreal}, \postcode{H3G1M8}, \state{QC}, \country{Canada}}}

\affil[3]{\orgdiv{School of Industrial and Systems Engineering}, \orgname{University of Oklahoma}, \orgaddress{%\street{Street}, 
 \city{Norman}, \postcode{73019}, \state{OK}, \country{United States}}}

%%==================================%%
%% Sample for unstructured abstract %%
%%==================================%%

\abstract{
Autonomous vehicles (AVs) can improve efficiency, reduce costs,
and enhance road safety. They optimize traffic flow, minimize congestion, and
support sustainability through shared mobility and reduced fuel consumption.
A key challenge in AV deployment is allocating vehicles to parking lots across
regions to meet fluctuating demand. Proper allocation reduces delays, lowers
costs, and boosts user satisfaction by ensuring timely vehicle availability. This
paper explores the impact of customer wait time patience on AV allocation models,
allowing prioritization of ride requests while balancing fleet efficiency and user
satisfaction. It also addresses the effectiveness of vehicle pooling in decentralized
service areas. We propose a mathematical model integrating vehicle distribution
and customer patience to maintain both efficiency and satisfaction. Results show
that adding more facilities initially reduces costs, but the benefits diminish with
more facilities. Increased customer patience improves inventory pooling benefits,
especially when no fixed costs are tied to facility operation.
}


%%\pacs[JEL Classification]{D8, H51}

%%\pacs[MSC Classification]{35A01, 65L10, 65L12, 65L20, 65L70}

\maketitle

\renewcommand{\thefootnote}{}%
\footnotetext{This is the Author’s Original Manuscript (AOM) of an article accepted for publication in \textit{International Journal of Logistics Systems and Management}, Inderscience Publishers. The article is currently listed among forthcoming articles at: \url{https://www.inderscience.com/info/ingeneral/forthcoming.php?jcode=ijlsm}.}

\section{Introduction}

Autonomous vehicles (AVs) are rapidly evolving technologies poised to transform the transportation industry. By eliminating the need for human drivers, these vehicles offer significant cost-saving opportunities for transportation companies as they strive to meet customer travel demands more efficiently \cite{bagloee2016autonomous, moradi2023last}. Furthermore, AVs reduce the risk of accidents caused by human error, as they operate without direct human intervention and rely on advanced sensing, decision-making, and control systems \cite{gharavi2024proactive, geisslinger2023maximum, sheikh2023collision, wang2020safety}. State-of-the-art AV systems employ vehicle-to-vehicle (V2V) and vehicle-to-infrastructure (V2I) communication to optimize traffic flow and minimize congestion, while adaptive routing algorithms enable real-time adjustments to enhance traffic efficiency \cite{stamadianos2023routing, guo2023sustainability}. These systems are already being piloted by companies like Waymo, Cruise, and Tesla for ride-sharing and logistics applications, demonstrating their feasibility in urban environments. AVs also support environmental sustainability by facilitating the adoption of shared mobility models and reducing fuel consumption \cite{bala2023review, pamidimukkala2023review, zhu2023comprehensive, singh2023adoption}.

The allocation of AVs to parking lots across different regions is critical to optimizing their operational efficiency and effectiveness. Appropriate allocation of AVs ensures that AVs are positioned where they are most likely to be demanded, minimizing service delays and enhancing user satisfaction. Proper distribution of AVs reduces empty vehicle travel, lowering energy consumption and operational costs and mitigating congestion in high-demand areas.

Effective allocation of AVs to parking lots across different regions considers various factors, including regional demand patterns, traffic conditions, and the availability of parking infrastructure \cite{zhang2023equilibrium, xie2023shared, sayarshad2023designing,choi2023analytical}. AV allocation optimization models can respond to fluctuations in demand, ensuring that vehicles are available when and where they are demanded on time. This is particularly important in urban areas, where demand for mobility services often varies by time of day and location. By addressing these variations, AV fleets can provide a more seamless and reliable transportation experience.

Incorporating customer wait time patience into AV allocation models is crucial for prioritizing ride requests and ensuring timely service delivery. Variations in customer tolerance levels enable the model to prioritize low-patience customers in critical situations. This consideration helps balance meeting immediate demand with maintaining fleet efficiency and user satisfaction.

The key contribution of this paper is the development of an advanced resource allocation model that integrates optimal distribution of AVs and customer wait time patience \cite{patel2023identifying, zhou2023understanding, ghandeharioun2023real}. By proposing an optimization model, this paper distributes AVs across parking lots and regions to ensure that AVs can respond effectively to customer demand with consideration of customer wait time patience. This enhances the adaptability of AVs, making them more robust while maintaining efficiency and customer satisfaction.

The remainder of the paper is as follows. Section \ref{background} gives an overview of the previous studies that support the assumptions of the optimization models for allocating AVs to customers. Section \ref{problem_definition} expresses the described system as an optimization model. Section \ref{proposed_heuristic_study} proposes a heuristic algorithm to solve the model. Section \ref{results} examines the performance of the proposed algorithm under various facility layout scenarios with different customer behaviors. Concluding remarks are offered in Section \ref{conclusion}.

\section{Literature review} \label{background} 

The routing and allocation of autonomous vehicles (AVs) to parking facilities plays a crucial role in enhancing the operational efficiency of transportation systems \cite{bahrami2020autonomous, zhao2023review, tang2024assessment, tian2023park, nourinejad2018designing, d2008future}. This is particularly vital in urban residential areas, where traffic congestion is common. Effective placement and movement of AVs are essential to maintaining the quality of transportation and delivery services \cite{guo2022managing, acquaviva2015customer, arslan2007autonomous, lam2016autonomous, abosuliman2021routing}. The challenge intensifies when customers exhibit varying levels of patience for waiting, highlighting the need for a model that effectively addresses transportation demands through AVs.

The deployment of AVs can be approached in various ways. One common strategy is to utilize centralized facilities, where vehicles are consolidated in one location to reduce the cost of opening multiple facilities \cite{ferrari2003car, vargas2008car}. While this strategy benefits from inventory pooling to minimize costs, it may lead to higher trip expenses for customers \cite{mason2003integrating, benjaafar2005benefits}. In fact, \cite{berlingerio2017graal} explored the impact of inventory pooling in transportation systems and proposed a data-driven approach to optimize carpooling by minimizing city-wide vehicle usage while enhancing passenger experience through social factors like interests and social connections. Eppen et al. \cite{eppen1979note} were pioneers in demonstrating the advantages of inventory pooling, showing that, under certain conditions, centralizing inventory can reduce costs by $\sqrt{n}$, where $n$ represents the number of locations. This cost-saving occurs because centralization smooths out demand fluctuations, with above-average and below-average demands canceling each other out.

However, customer wait time tolerance adds complexity to AV routing and allocation. While some customers are patient and willing to wait for a vehicle from a distant location, others expect prompt service. In cases where customers exhibit patience, fewer facilities can be open, leveraging inventory pooling to reduce vehicle stationing costs. Conversely, for customers with low tolerance for waiting, spreading vehicles across a broader area may be more effective in ensuring prompt service \cite{strohle2019leveraging}. Moreover, solving vehicle assignment problems is challenging due to their intricate nature. Several studies have applied heuristic approaches, such as hybrid algorithms \cite{toffolo2019heuristics}, to improve efficiency in addressing these problems.

The effectiveness of inventory pooling is influenced by demand patterns. Eppen's model demonstrated that inventory pooling is advantageous only under specific demand conditions, such as when demand is negatively correlated across locations. In contrast, when demand is positively correlated, pooling offers no benefit. Furthermore, when demand follows a heavy-tailed distribution, pooling becomes less effective. Other factors, such as service levels, process variability, and control policies, also influence the success of inventory pooling \cite{vargas2008car, ferrari2003car, benjaafar2005benefits, swinney2012inventory, bimpikis2016inventory}. Previous research has explored the tradeoff between inventory pooling and customer proximity, with the facility location problem providing a framework for finding the optimal number and location of facilities to minimize costs \cite{cournuejols1990uncapacitated}. In the context of AVs, this problem becomes more complex as the chance of losing a sale due to distance must also be considered. This paper introduces a new factor: the impact of losing a sale as the distance between AVs and customers increases.

The facility location problem has multiple variations, such as the capacitated facility location problem, where facilities have limits on how much inventory they can hold. This problem is particularly relevant to AV allocation, as there is a fixed number of vehicles available. Benjaafar et al. \cite{benjaafar2017inventory} studied inventory repositioning in on-demand product rentals, examining the tradeoff between repositioning costs and the risk of losing sales due to poorly located inventory. Their model allowed customers to rent a product without reservation and return it to any eligible location. This contrasts with our model, where the destination and trip duration are known in advance, focusing on vehicle allocation and customer demand.

Braverman et al. \cite{braverman2019empty} investigated empty vehicle routing in systems similar to ride-sharing networks. This research is relevant to our work, as it examines the strategy of deciding whether a vehicle should wait at its current location or move to another area. In their model, a decision matrix determines the probability of a vehicle moving between locations after completing a trip. In our study, this decision is crucial, as the distance to customers depends on the location of the vehicles. For our simulation study, vehicles always return to their home facility after completing a trip. This decision ensures that the primary tradeoff between pooling and proximity remains the focus, as the vehicles are initially distributed across various locations.

In this paper, we focus on determining the optimal spatial distribution of autonomous vehicles, considering both customer demand and the unique characteristics of AVs. Unlike ride-sharing, which involves human drivers, our model removes the human factor, enabling a clearer examination of the trade-offs between pooling and proximity. We also incorporate customer patience levels, recognizing that different customers have varying tolerance for wait times. While previous research has explored inventory repositioning in on-demand product rentals \cite{benjaafar2017inventory}, our model is distinct in that the destination and trip duration are known in advance, allowing us to specifically focus on customer demand, vehicle allocation, and the impact of customer patience.

\section{Problem definition and formulation}\label{problem_definition}

This section explains the model's primary components, including sets and parameters, decision variables, constraints, and the objective function. In the proposed model, we introduce a grid of locations. These locations host a group of customers traveling to various destinations denoted across different time units. The proposed model incorporates different facilities for stationing AVs. The decision-making process involves determining whether a facility should remain closed or open. Furthermore, the model addresses the allocation of available AVs to fulfill customer requests, aiming to minimize the total system cost. This includes variable expenses such as distance travel costs, lost sale costs, and fixed costs like keeping facilities open. In this model, lost sale costs may occur due to vehicle shortages at specific times for particular customers or customer patience constraints.
The model's notations are outlined in Table \ref{tabel1}. 

First, we have the set of time units ($\mathcal{T}$), representing the discrete time units over which decisions must be made. Each time unit corresponds to a specific moment when vehicle assignments, relocations, and facility operations are updated. The set of customers ($\mathcal{M}$) includes all individuals who may request transportation services, with each customer indexed by $m \in \mathcal{M}$. The set of locations ($\mathcal{V}$) represents the geographic areas where customers and AVs are situated, and each location is indexed by $v \in \mathcal{V}$. There are subsets of available customers and vehicles at each location and time. The set of available customers at location $v$ and time $t$ ($\mathcal{M}_{vt}$) represents customers who are present at a specific location and time, awaiting service. Similarly, the set of available AVs at location $v$ and time $t$ ($\mathcal{L}_{v,t}$) indicates the number of vehicles ready for deployment at a given location and time.

The facility capacity ($F_t$) determines the overall capacity to keep certain facilities open at time $t$, allowing the storage and maintenance of AVs. Additionally, a large constant value, denoted as $N$, is used in constraints to enforce certain decisions. Distance between locations ($d_{uv}$) represents the physical distance between any two locations $u$ and $v$, and travel time between locations ($t_{uv}$) defines the time it takes for a vehicle to travel between those exact locations. Both distance and travel time are essential for calculating transportation costs and ensuring timely service. In terms of costs, the lost sale cost per unit at time $t$ ($\mathcal{K}_{1,t}$) refers to the penalty for not serving a customer due to vehicle unavailability. The distance travel cost per unit at time $t$ ($\mathcal{K}_{2,t}$) captures the cost incurred for vehicle movement between locations. Additionally, the fixed cost for keeping a facility open at time \( t \) (\(\mathcal{K}_{3,t}\)) represents the cost of operating a vehicle facility during a specific period. These time-dependent costs account for variations in operational factors such as staffing, maintenance, and demand fluctuations, which change over time. Allowing facilities to alternate between open and closed is realistic in scenarios with fluctuating demand. This flexibility helps reduce costs by shutting down facilities during off-peak times, reflecting operational efficiency in dynamic environments like transportation or logistics, where continuous facility operation may not always be necessary.
Another critical factor is the patience radius for customer $m$ at location $v$ and time $t$ ($S_{m,v,t}$), which describes the maximum distance a customer is willing to wait for a vehicle. If no vehicle arrives within this radius, the demand is considered unmet.

Several decision variables dictate how the system operates. The binary variable $x_{v,m,k}^t$ takes the value of 1 if a customer $m$ located at $v$ at time $t$ is assigned to an autonomous vehicle at the exact location to be dropped off at location $k$. Otherwise, this variable is set to 0. Another binary variable, $y_{mv,cu,k}^t$, is one if a customer $m$ at location $v$ at time $t$ is assigned to a vehicle $c$ at a different location $u$, and the customer is dropped off at location $k$; it is zero otherwise. For vehicle relocations, the binary variable $z_{c,u,v}^t$ takes the value of 1 if an autonomous vehicle $c$ located at $u$ is returned or relocated to another location $v$ at time $t$. If no relocation occurs, it remains 0. The variable $n_{v,t}$ represents the number of available AVs at location $v$ at time $t$, capturing the system’s ability to meet demand at specific points. Also, the binary variable $w_{m,v,t}$ is set to 1 if a lost sale occurs for customer $m$ at location $v$ at time $t$, indicating that the demand could not be met due to a lack of available vehicles. Finally, $f_{v,t}$ is a binary variable equal to 1 if a facility is opened at location $v$ to store AVs at time $t$, and zero if the facility remains closed. The mathematical model for the addressed problem is presented as follows:
 
% The simulation section considers the impact of inventory pooling and the characteristics of the demand distribution.
% Table 1 describes the details of the model. The simulation section considers the effects of Inventory pooling and the nature of the demand distribution. 

% \subsection{Model Description}
\begin{table}[ht!]
\centering
\caption{Notations and variables of the proposed model}\label{tabel1}
%\begin{adjustbox}{width=1\textwidth}
\begin{tabular}{lp{11cm}}
\hline
\textbf{Sets} & \textbf{Description} \\
%\midrule
%Data & \\
\hline
$\mathcal{T}=\{1,2,3,..T\}$ & Set of time units\\
$\mathcal{M}=\{1,2,3,..M\}$ & Set of customers\\
% $i$                & x coordinate of location  \\
% $j$                & y coordinate of location  \\
$\mathcal{V}$ & Set of locations \\
$\mathcal{M}_{vt}$ & Set of available customers at location $v\in \mathcal{V}$ at time $t\in \mathcal{T}$ \\
$\mathcal{L}_{v,t}$ & Set of available  AVs at location $v\in \mathcal{V}$ at time $t\in \mathcal{T}$ \\
% $\mathcal{D}_(1,t)$ & Set of customers who can be picked by vehicles at their locations\\ & at time $t\in \mathcal{T}$\\
% $\mathcal{D}_(2,t)$ & Set of customers who can be picked by vehicles at other locations rather \\ & than their location at time $t\in \mathcal{T}$\\
\hline
\textbf{Parameters} & \\
\hline
$F_t$&	The overall capacity to maintain certain facilities open at time $t \in \mathcal{T}$\\
$N$ & A large constant value\\
% $v$                & Location of customer or vehicle as a two dimensional vector $(x,y)$\\
% $c_{v,t}$ & The index of existing autonomous vehicle at location $v \in \mathcal{V}$ at time $t \in \mathcal{T}$ \\
$d_{uv}$ & Distance between the locations $u,v \in\mathcal{V}$ \\
$t_{uv}$ & Travel time between the locations $u,v \in\mathcal{V}$ \\
$\mathcal{K}_{1,t}$ & Lost sale cost per unit at time $t\in \mathcal{T}$\\
$\mathcal{K}_{2,t}$ & Distance travel cost per unit at time $t\in \mathcal{T}$ \\
$\mathcal{K}_{3,t}$ & Fixed cost for keeping a facility open at time $t\in \mathcal{T}$\\
$S_{m,v,t}$ & Patience radius for customer $m\in \mathcal{M}$ at location $v\in \mathcal{V}$ at time $t\in \mathcal{T}$\\
\hline
\textbf{Decision Variables} & \\
\hline
$x_{v,m,k}^t$ & $1$ if customer $m\in \mathcal{M}_{v,t}$ at location $v\in\mathcal{V}$ at time $t \in \mathcal{T}$ is assigned to an autonomous vehicle at their exact location to be dropped off at location $k\in\mathcal{V}$; $0$ otherwise.\\
$y_{mv,cu,k}^t$ & $1$ if customer $m\in \mathcal{M}_{v,t}$ at location $v\in\mathcal{V}$ at time $t\in \mathcal{T}$ is assigned to the autonomous vehicle $c \in \mathcal{L}_{u,t}$ at location $u \in \mathcal{V}$ to be dropped off at location $k\in\mathcal{V}$; $0$ otherwise.\\
$z_{c,u,v}^t$ & $1$ if the autonomous vehicle $c \in \mathcal{L}_{u,t}$ at location $u \in \mathcal{V}$ is returned to a different location $v\in \mathcal{V}$ at time $t\in \mathcal{T}$; 0 otherwise.\\
$n_{v,t}$ & Number of available AVs at location $v\in \mathcal{V}$ at time $t\in \mathcal{T}$\\
$w_{m,v,t}$ & $1$ if there is a lost sale unit for customer $m\in \mathcal{M}$ at location $v\in \mathcal{V}$ at time $t\in \mathcal{T}$; $0$ otherwise. \\
$f_{v,t}$ &	$1$ if a facility is opened to store AVs at location $v\in \mathcal{V}$ at time $t\in \mathcal{T}$; $0$ otherwise.\\
\hline
\end{tabular}
%\end{adjustbox}
\end{table}

\begin{flalign}\label{eq:objectivefunction}
\notag
& \underset{x,y,z,w,f,n}{\mathrm{min.}} \quad \sum_{v \in V} \sum_{t \in \mathcal{T}} \mathcal{K}_{1,t}\sum_{m \in \mathcal{M}_{v,t}} w_{m,v,t} + \sum_{v,k \in V: k \neq v}\sum_{t \in \mathcal{T}} \mathcal{K}_{2,t}\sum_{m \in \mathcal{M}_{v,t}} x_{v,m,k}^t(d_{vk} + d_{kv}) \\ \notag
&+\sum_{v,u,k \in V: v \neq u, k \neq u, k \neq v} \sum_{t \in \mathcal{T}}  \mathcal{K}_{2,t}\sum_{c \in \mathcal{L}_{u,t}} \sum_{m \in M_{v,t}} y_{mv,cu,k}^t(d_{uv}+d_{vk}+d_{ku}) \\ 
&+ \sum_{v \in V}\sum_{t \in \mathcal{T}} \mathcal{K}_{3,t}{f_{v,t}} 
\end{flalign}
\begin{flalign}\label{con1}
\notag
& \text{s.t.} \\ 
& \sum_{m \in \mathcal{M}_{u,t}}\sum_{k \in \mathcal{V}: k \neq u} x_{u,m,k}^t +\sum_{v,k \in \mathcal{V}: v\neq k , v \neq u, k \neq u}\sum_{m \in \mathcal{M}_{v,t}}\sum_{c \in \mathcal{L}_{u,t}} y_{mv,cu,k}^t\leq n_{u,t}f_{u,t},\forall u \in \mathcal{V},\forall t \in \mathcal{T}
\end{flalign}
\begin{flalign}\label{con2}
\notag &n_{v,t}=n_{v,t-1}-\sum_{m \in \mathcal{M}_{v,t}}\sum_{k \in \mathcal{V}: k \neq v} x_{v,m,k}^t-\sum_{u,k \in V: u \neq k,k \neq v, u \neq v}\sum_{m \in \mathcal{M}_{v,t}}\sum_{c \in \mathcal{L}_{u,t}} y_{mv,cu,k}^t\\&+\sum_{c \in \mathcal{L}_{u,t}} z_{c,u,v}^t,\forall v \in \mathcal{V}, \forall t \in \{2,3,..T\}
\end{flalign}
\begin{flalign}\label{con3}
& \sum_{v,k \in \mathcal{V}: k \neq v, k \neq u,v \neq u} \sum_{m \in \mathcal{M}_{v,t}}y_{mv,cu,k}^t\leq z_{c,k,u}^{(t+t_{uv}+t_{vk}+t_{ku})} \quad \forall c \in \mathcal{L}_{u,t},\forall u \in \mathcal{V},\forall t \in \mathcal{T}
\end{flalign}
\begin{flalign}\label{con4}
&\sum_{u,k \in \mathcal{V} : u \neq k, u \neq v }  \sum_{c \in \mathcal{L}_{u,t}} y_{mv,cu,k}^t\leq 1 \quad \forall v \in \mathcal{V},\forall m \in \mathcal{M}_{v,t}, \forall t \in \mathcal{T}
\end{flalign}
\begin{flalign}\label{con5}
& \sum_{v,k \in \mathcal{V}: v \neq k , v \neq u,k \neq u} \sum_{m \in \mathcal{M}_{v,t}} y_{mv,cu,k}^t\leq 1 \quad \forall u \in \mathcal{V},\forall c\in \mathcal{L}_{u,t}, \forall t \in \mathcal{T}
\end{flalign}
\begin{flalign}\label{con6}\notag
& N(1- \sum_{u \in V :u \neq v, u \neq k} \sum_{c \in \mathcal{L}_{u,t}}y_{mv,cu,k}^t)+ \sum_{u \in V :u \neq v, u \neq k } \sum_{c \in \mathcal{L}_{u,t}}y_{mv,cu,k}^td_{uv}\leq S_{m,v,t}(1-w_{m,v,t})\\&+N w_{m,v,t}, \forall m \in \mathcal{M}_{v,t},\forall v,k \in \mathcal{V},\forall t \in \mathcal{T};  v \neq k
\end{flalign}
\begin{flalign}\label{con7}
&\sum_{t \in \mathcal{T}}x_{v,m,k}^t+\sum_{t \in \mathcal{T}}w_{m,v,t}+\sum_{u \in \mathcal{V}: u \neq v, u \neq k}\sum_{c \in \mathcal{L}_{u,t}}\sum_{t \in \mathcal{T}} y_{mv,cu,k}^t=1 \quad\forall v,k \in \mathcal{V},\forall m \in \mathcal{M}_{v,t}; v \neq k
\end{flalign}
\begin{flalign}\label{con8}
% &\sum_{v,k \in \mathcal{V}: v\neq k , v \neq u, k \neq u}\sum_{m \in \mathcal{M}_{v,t}}\sum_{c \in \mathcal{L}_{u,t}} y_{mv,cu,k}^t\leq n_{u,t}f_{u,t}\quad\forall u \in \mathcal{V},\forall t \in \mathcal{T}, \\ \label{con9}
& \sum_{v \in V}{f_{v,t}} \leq F_t \quad \forall t \in \mathcal{T}
\end{flalign}
\begin{flalign}\label{con10}
& x_{v,m,k}^t\in \{0,1\} \quad\forall v,k \in \mathcal{V},\forall m \in \mathcal{M}_{v,t},\forall t \in \mathcal{T}
\end{flalign}
\begin{flalign}\label{con11}
& y_{mv,cu,k}^t\in \{0,1\}\quad  \forall v,u,k \in \mathcal{V},\forall m \in \mathcal{M}_{v,t}, \forall t \in \mathcal{T},\forall c \in \mathcal{L}_{v,t},\forall t \in \mathcal{T}; v\neq k , v \neq u, k \neq u
\end{flalign}
\begin{flalign}\label{con12}
& z_{c,u,v}^t \in \{0,1\}\quad \forall u,v \in \mathcal{V}, \forall c \in \mathcal{L}_{v,t},\forall t \in \mathcal{T}; u\neq v
\end{flalign}
\begin{flalign}\label{con13}
&w_{m,v,t} \in \{0,1\}\quad\forall v \in \mathcal{V},\forall m \in \mathcal{M}_{v,t}, \forall t \in \mathcal{T}
\end{flalign}
\begin{flalign}\label{con14}
&f_{v,t}\in \{0,1\} \quad \forall v \in \mathcal{V} , \forall t \in \mathcal{T}
\end{flalign}
\begin{flalign}\label{con15}
&n_{v,t}\in Z^+ \quad \forall v \in  \mathcal{V}, \forall t \in \mathcal{T}
\end{flalign}

The objective function \eqref{eq:objectivefunction} minimizes the overall system cost, including the total lost sale costs, total distance travel costs (either by a vehicle in the exact customer location to the destination or from a different location), and the total fixed cost of keeping facilities open.
Constraint \eqref{con1} assigns the available AVs at each location and time unit to the customers at the same or different locations. This constraint guarantees that the total number of AVs assigned to customers at different locations does not surpass the total number of available AVs at that location for each time unit.
Constraint \eqref{con2} determines the total number of available AVs at each time unit for each location. This calculation entails subtracting the total number of vehicles leaving that location for customer pick-up, whether at the exact or different locations, from the total available AVs at the previous time unit. It also accounts for the vehicles returning to this location at the same time unit. Constraint \eqref{con3} guarantees that if a specific autonomous vehicle at a particular location and time has been assigned to a customer at a different location, it must be returned to its initial location after completing its trip.  This is to simplify fleet management, reduce repositioning costs, and ensure centralized maintenance. Constraint \eqref{con4} specifies that each customer can be transported by at most one autonomous vehicle during each time unit. Constraint \eqref{con5} states that each autonomous vehicle can be assigned to at most one customer during each time unit. Ridesharing is not considered as vehicle capacity is assumed to meet the demand of a single customer, and allowing multiple customers would add unnecessary complexity to the model. Constraint \eqref{con6} represents that when a request is received from a customer, an autonomous vehicle can be assigned. If the customer is patient, they can wait for the assigned autonomous vehicle to pick them up. A lost sale cost is incurred if the customer is impatient and unwilling to wait. Constraint \eqref{con7} specifies that each demand must be fulfilled by an autonomous vehicle existing precisely at the exact location, or it can be met after some time by assigning an available vehicle from a different location, or it is lost due to customer impatience. Constraint \eqref{con8} also states that the number of open facilities for dispatching vehicles to different customers at each time unit should not surpass the total allowed capacity. This reflects practical limitations such as staffing, energy consumption, and maintenance costs, ensuring resources are utilized efficiently. It also accounts for fluctuating demand, avoiding the unnecessary expense of keeping underutilized facilities open, and promotes cost-effective operations in dynamic service environments. Constraints \eqref{con10}-\eqref{con15} are the domain constraints.

\section{Solution method} \label{proposed_heuristic_study}

The proposed model (\ref{eq:objectivefunction}-\ref{con15}) in Section \ref{problem_definition} is computationally complex, primarily due to its categorization as a non-linear assignment problem. Non-linear assignment problems, including the generalized assignment problem, are well-documented in the literature as NP-hard~\citep{fisher1986multiplier}. These problems are characterized by their combinatorial nature, which requires determining optimal assignments between two sets (e.g., vehicles and customers) under given constraints. The difficulty lies in the exponential growth of the solution space with the problem size, making it computationally intractable to solve large-scale instances using exact optimization methods. 

Due to the limitations of exact methods for large and dynamic systems, heuristic and metaheuristic approaches have gained significant attention. Heuristic methods are particularly effective in producing high-quality solutions within a reasonable computational time, making them practical for real-time applications. Inspired by this motivation, we propose a heuristic method named \textit{``A Dynamic Matching Algorithm for Daily Reset Autonomous Vehicle Allocation''} to approximate solutions for the problem. The method leverages principles of dynamic decision-making and prioritization to balance computational efficiency and solution quality.

The dynamic matching algorithm is structured around real-time evaluation of customers and vehicles in a service system. At its core, the heuristic evaluates customers at each location based on their ``patience'', which is defined as the maximum amount of time they are willing to wait for service. This parameter serves as the primary driver for prioritization in the allocation process. Customers with lower patience are given higher priority, reflecting their increased likelihood of abandoning the system if not served promptly. 

At each discrete time unit, the heuristic performs a sequence of operations. It begins by checking the availability of autonomous vehicles at the customer's current location. If a vehicle is available, it is immediately assigned to the customer, ensuring that no other customer is being served at the same location simultaneously. This constraint guarantees that each vehicle is dedicated to a single task at any given time. If no vehicle is available at the exact location, the algorithm performs a broader search across other locations to identify the nearest available vehicle. This search process considers the distance between the customer and potential vehicles, as well as the customer's patience and waiting radius. The waiting radius represents the maximum distance a customer is willing to tolerate for a vehicle to arrive. If the nearest available vehicle is beyond this threshold, the customer cancels their request, resulting in a lost sale. Conversely, if the vehicle is within the waiting radius, it is dispatched to the customer. 

This matching mechanism introduces a trade-off between customer satisfaction and operational efficiency. By prioritizing customers with lower patience and ensuring vehicles are assigned based on proximity, the algorithm strives to minimize both customer churn and vehicle idle times. Additionally, the heuristic accounts for real-world constraints such as vehicle travel times, pickup and drop-off operations, and system-wide resource availability.

A critical aspect of the proposed heuristic is its consideration of vehicle movement and availability. Each autonomous vehicle is assumed to travel at a fixed speed, set to one block per minute in this study. Once a vehicle is matched with a customer, it begins traveling to the customer's location. Upon arrival, the vehicle picks up the customer and transports them to their desired destination. After completing the trip, the vehicle returns to its facility, where it becomes available for future assignments. Importantly, a vehicle remains unavailable during the entire cycle of pickup, transportation, and return, ensuring that it cannot be reassigned until it completes its current task.

To reflect operational realities, the algorithm prohibits vehicles from serving multiple customers simultaneously or abandoning an ongoing task to serve a closer customer. This constraint ensures that each vehicle's commitment to a customer is honored. Furthermore, customers are modeled with individual patience levels and waiting radii to capture the variability in service expectations and tolerances.

One of the defining features of the heuristic is its daily reset mechanism. At the end of each operational day, vehicles with ongoing tasks complete their assignments and return to their respective facilities. Vehicles that remain idle throughout the day also return to their facilities. This reset provides an opportunity to reinitialize the system, allowing for a balanced redistribution of resources and preparation for the next day's operations. The reset ensures fairness and efficiency, as vehicles are evenly distributed across facilities and ready to handle new requests. The daily reset also simplifies the modeling of vehicle availability and resource allocation, as the system starts each day with all vehicles in known locations. This feature is particularly advantageous for scenarios with fluctuating demand patterns, as it prevents resource imbalances from persisting across multiple days. Additionally, it provides a natural point for system performance evaluation and adjustment, enabling operators to fine-tune parameters such as vehicle deployment strategies, customer prioritization criteria, and waiting radius thresholds.

The proposed heuristic offers several advantages over traditional optimization methods. First, its dynamic nature allows it to adapt to real-time changes in customer demand and vehicle availability. By prioritizing customers based on patience levels and proximity to available vehicles, the algorithm ensures a high level of responsiveness and service quality. Second, its computational efficiency makes it suitable for large-scale applications, where exact methods would be prohibitively time-consuming. Third, the incorporation of practical constraints, such as vehicle travel times and commitment policies, enhances its applicability to real-world autonomous vehicle systems. Potential applications of the proposed heuristic extend beyond autonomous vehicle allocation. Similar approaches could be adapted for ride-sharing services, delivery logistics, and emergency response systems, where dynamic resource allocation and real-time decision-making are critical. The daily reset mechanism, in particular, is well-suited for systems with cyclical demand patterns, such as urban transit networks and shared mobility platforms.

In summary, the proposed dynamic matching algorithm provides a robust and scalable solution to the problem of autonomous vehicle allocation under time and space constraints. Its design reflects the complexities of real-world operations while balancing computational efficiency and solution quality. The combination of dynamic prioritization, distance-based matching, and a daily reset mechanism makes it a versatile tool for addressing challenges in dynamic service systems.

\begin{algorithm}
%\footnotesize
\caption{A Dynamic Matching Algorithm for Daily Reset Autonomous Vehicle Allocation}\label{Algorithm1}
\begin{algorithmic}[1]
% \State Set time of the system to $\boldsymbol{0}$
\State Initiate $\boldsymbol{n}$ scenarios with $\boldsymbol{m}$ open facility 
at random locations in each scenario
\State Distribute the total $\boldsymbol{l}$ available AVs randomly in each location 
\For {Every facility layout Scenario}
\For{Every open facility location at each facility layout scenario at each time $\boldsymbol{t}$}
\State Give the available customers priority based 
\State on their patience level (less patient customers take higher priority)
\For{Every customer at each level of priority}
    \State Look for an available vehicle at the customer's present location
    \If{Found an available vehicle and there is no other prior customer at that location }
        \State Match the vehicle with the customer and make the status
        \State of vehicle \textit{``traveling to the customer''} 
        \State Add the customer to set of \textit{``satisfied demands''} 
    \Else
        \State Look in other locations for the closest potential vehicle
        \If{The closest vehicle is within the patient radius and there are no other customers 
        \State with higher priority in the system at that time}
            \State Move the vehicle to pick up the customer and make the status
        \State of vehicle \textit{``traveling to the customer''} 
        \Else
        % \State Increase time of the system by $\boldsymbol{1}$
        \State Reject the sale, a lost sale cost occurs
         \State Add the customer to the set of ``rejected demands''.
        \EndIf
    \EndIf
    \EndFor

    % \State Each vehicle can move 1 block per minute
\For {Every assigned vehicle in the current open facility}
     \State Check the status \textit{(``idle'', ``traveling to the customer'', ``returning to the facility'')} 
     \State of the assigned vehicles at each time 
      \If{ A vehicle finished its trip at that time}
      \State Change the vehicle status to \textit{``idle''} at the facility 
    \ElsIf{A vehicle started its trip to a customer at that time}
      \State Change the vehicle status to \textit{``traveling to the customer''} at the facility 
    \ElsIf{A vehicle is idle at each location at that time}
      \State Maintain its status
    \EndIf  
\EndFor
\State Update the total cost of the system 
\EndFor
\State Reset the system at the end of the day 
\EndFor
\State Find the scenario with the minimum total cost 
\end{algorithmic}
\end{algorithm}

\section{Computational results}\label{results}

The proposed heuristic was coded in Java and compiled on a personal laptop using the software application Geany. We conducted multiple case scenarios and solved them using our suggested heuristic. In each scenario, we assume that the facility locations and total number of available vehicles at each location are known. To ensure consistency across scenarios, we utilized a uniform $7\times7$ two-dimensional grid for generating the locations of vehicles and customers.
We also maintain consistency in model parameters, namely $\mathcal{K}_{1,t}$, $\mathcal{K}_{2,t}$, and $\mathcal{K}_{3,t}$, $|\mathcal{M}|$, across all scenarios. In our experiments, each scenario spanned 100 days, and each day is assumed to be 600 units of time.  We assume 100 customers in the system, and the presence of a customer at each block is modeled using a Bernoulli distribution with a success probability of $5\%$.

\subsection{Facility Layout Configurations}
\begin{figure}[ht!]
\centering
\begin{tabular}{ccc}
\includegraphics[width=0.3\textwidth]{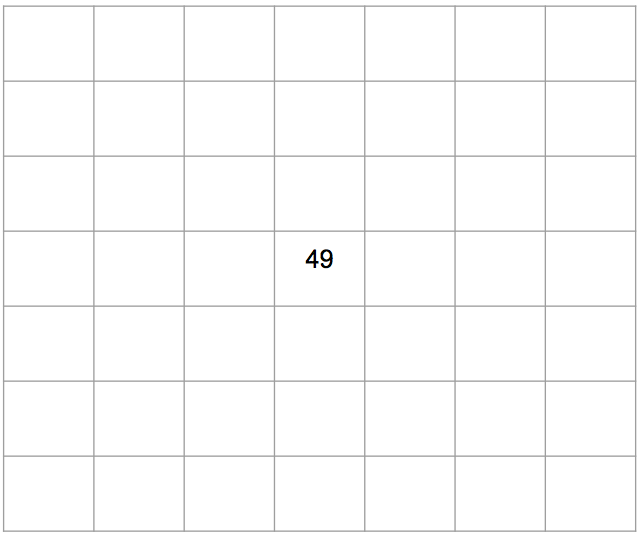}
& \includegraphics[width=0.3\textwidth]{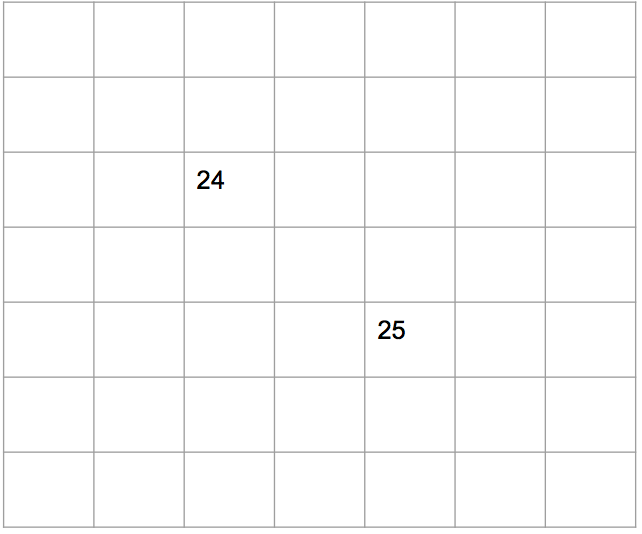}
& \includegraphics[width=0.3\textwidth]{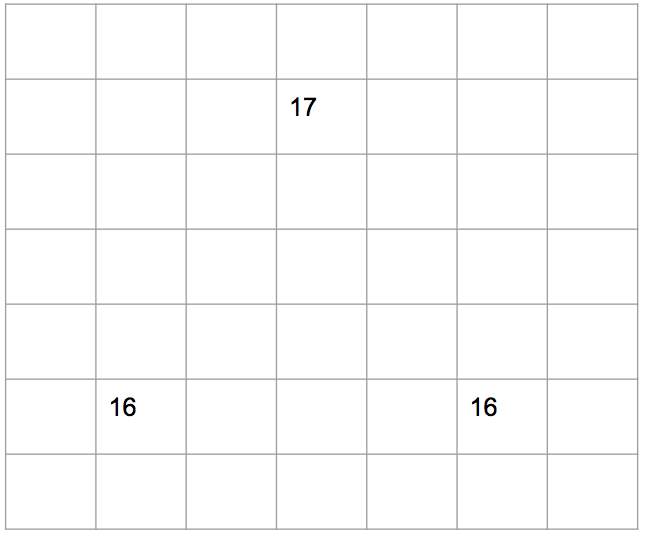} \\
Layout 1
& Layout 2
& Layout 3 \\
& & \\
\includegraphics[width=0.3\textwidth]{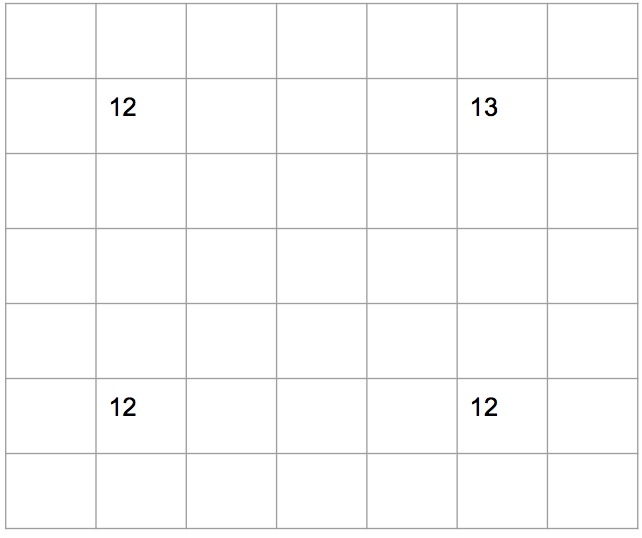} &
\includegraphics[width=0.3\textwidth]{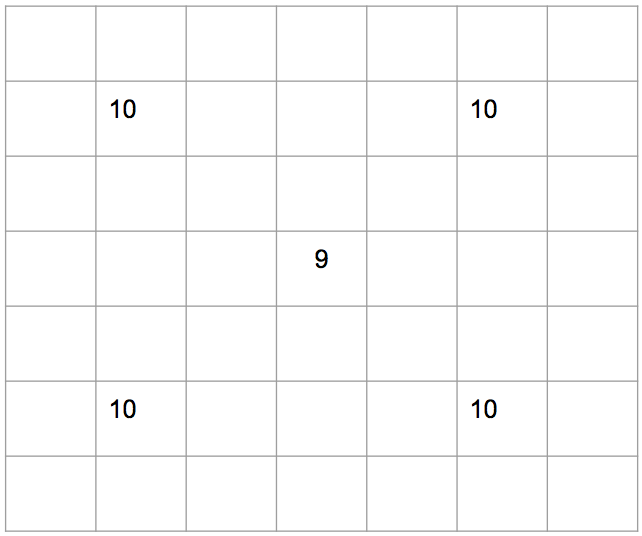} &
\includegraphics[width=0.3\textwidth]{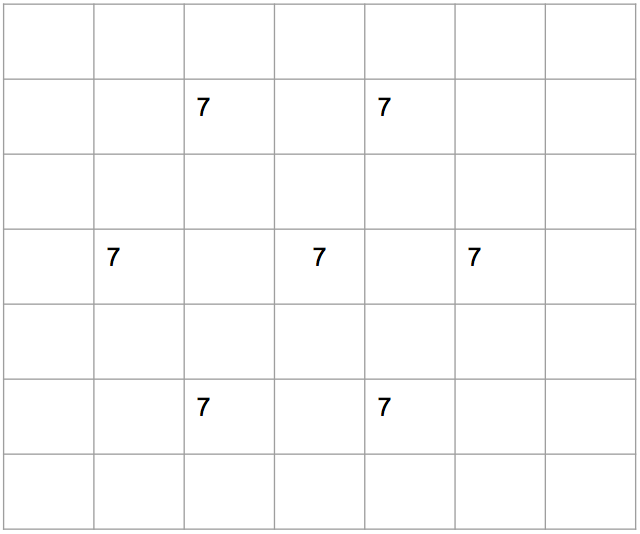}\\
Layout 4 &
Layout 5 &
Layout 6 \\
& & \\
\includegraphics[width=0.3\textwidth]{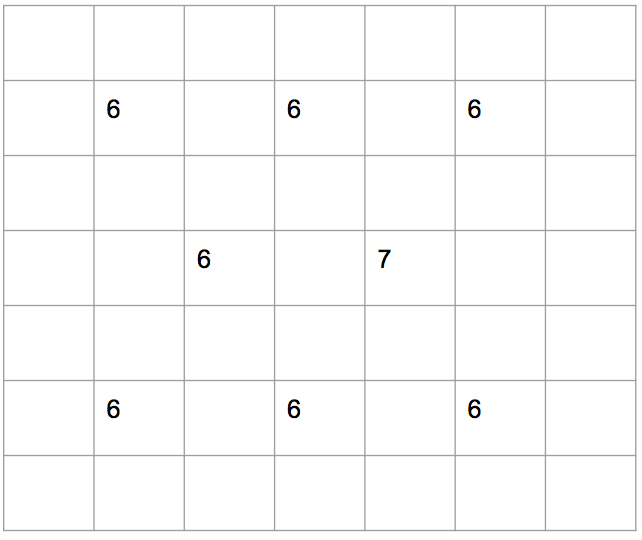} &
\includegraphics[width=0.3\textwidth]{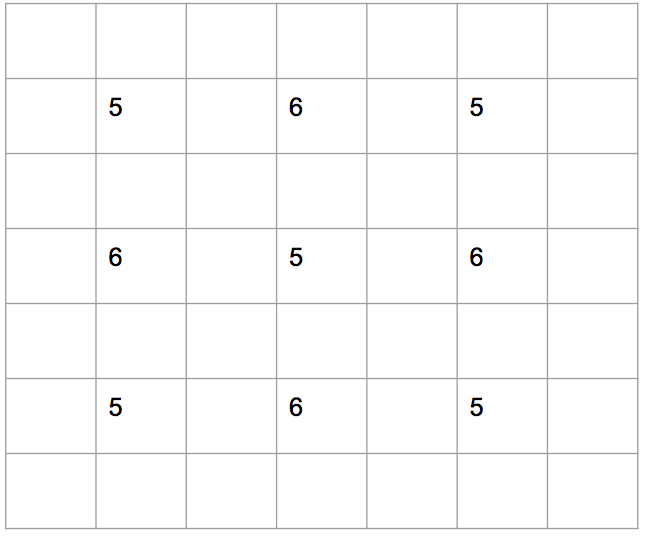} &
\includegraphics[width=0.3\textwidth]{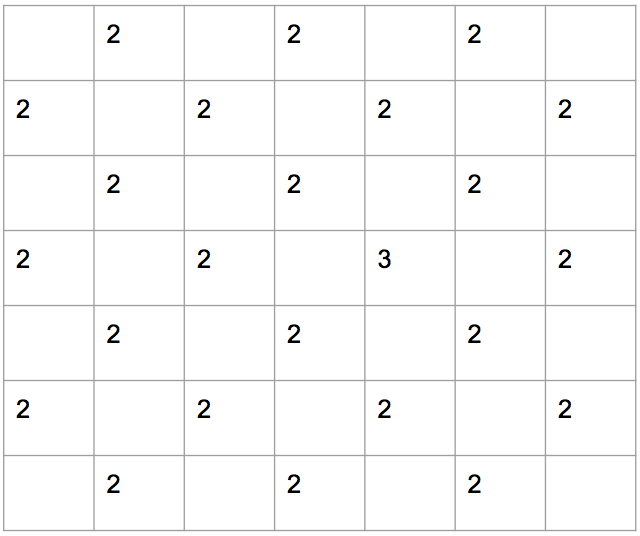} \\
Layout 7 &
Layout 8 &
Layout 9 \\
& & \\
\includegraphics[width=0.3\textwidth]{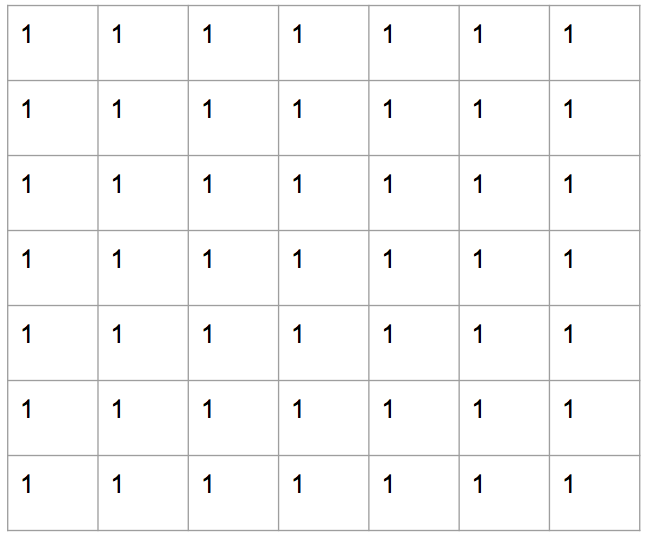} & & \\
Layout 10 & & \\
\end{tabular}
\caption{Different facility layouts}
\label{figconfigs}
\end{figure}

\begin{figure}[ht!]
\centering
\begin{tabular}{ccc}
\includegraphics[width=0.3\textwidth]{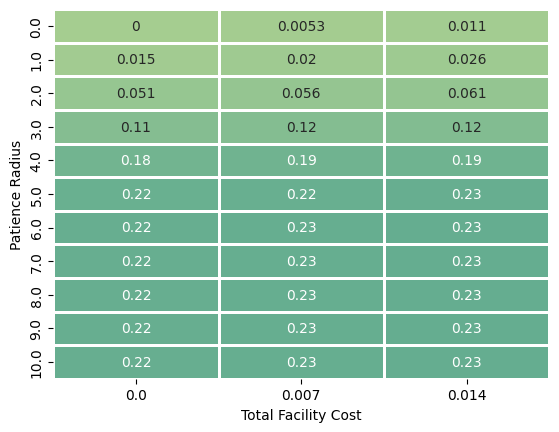}
& \includegraphics[width=0.3\textwidth]{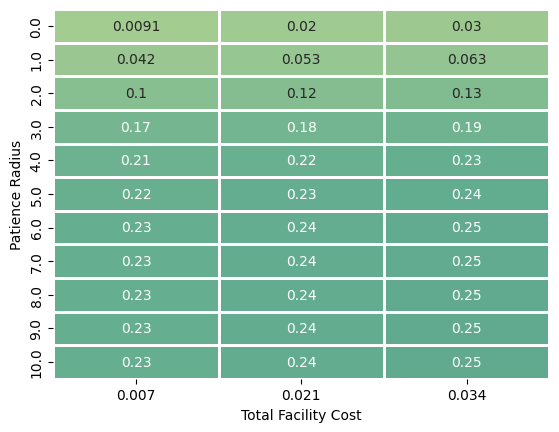}
&\includegraphics[width=0.3\textwidth]{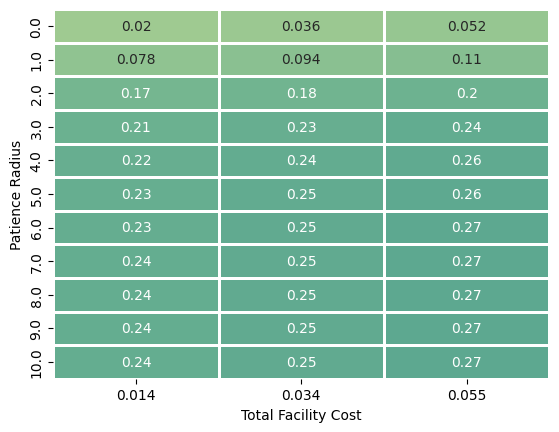} \\
Layout 1 &
Layout 2 &
Layout 3 \\
\includegraphics[width=0.3\textwidth]{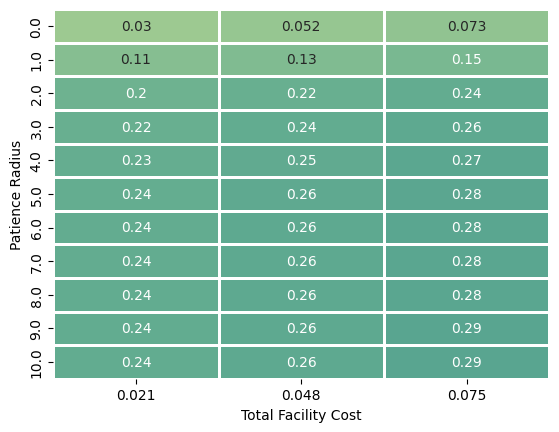} &
\includegraphics[width=0.3\textwidth]{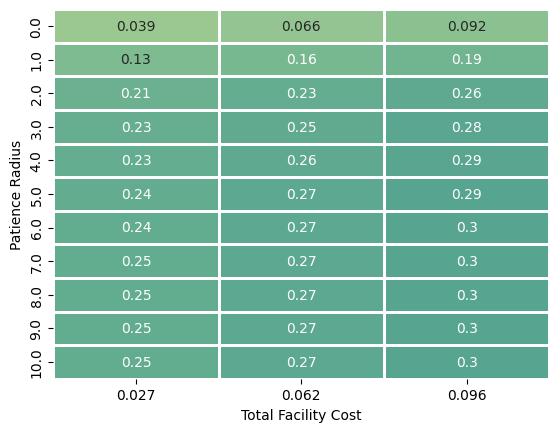} &
\includegraphics[width=0.3\textwidth]{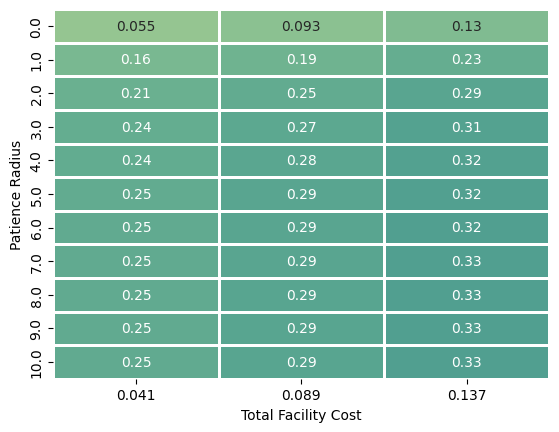}\\
Layout 4 &
Layout 5 &
Layout 6 \\
& & \\
\includegraphics[width=0.3\textwidth]{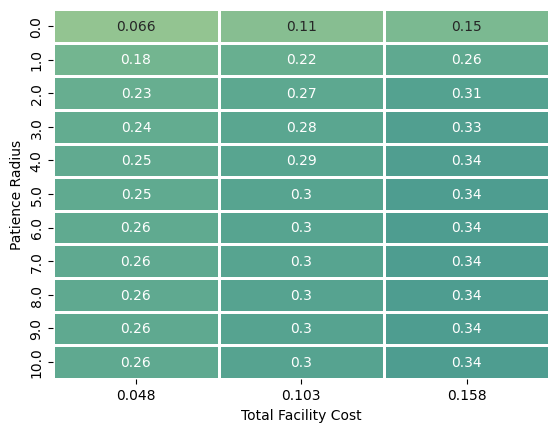} &
\includegraphics[width=0.3\textwidth]{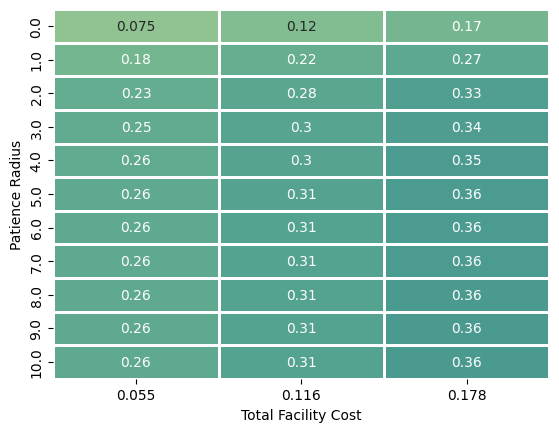} &
\includegraphics[width=0.3\textwidth]{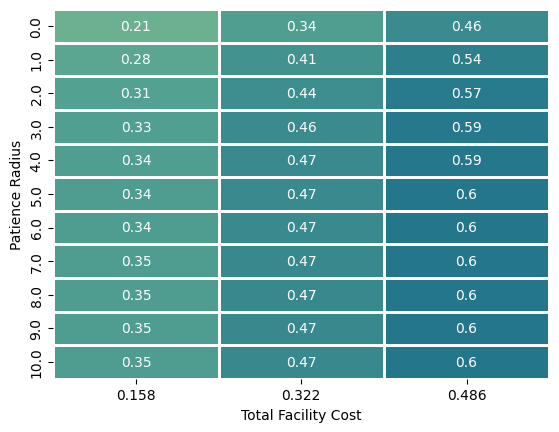} \\
Layout 7 &
Layout 8 &
Layout 9 \\
& & \\
\includegraphics[width=0.3\textwidth]{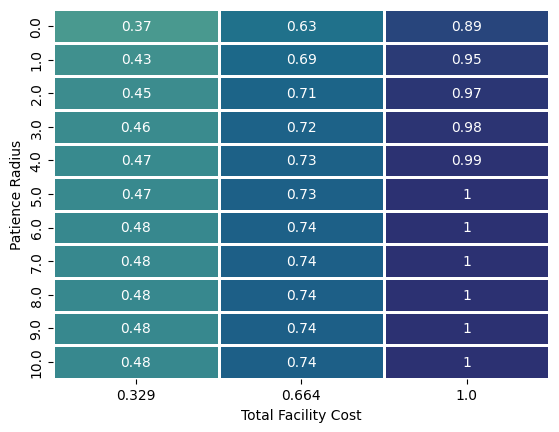} & & \\
Layout 10 & & \\
\end{tabular}
\caption{Normalized total cost of the system for different facility layouts}
\label{figconfigstotalcost}
\end{figure}

We designed different configurations to assess the impact of various facility opening strategies on the system's total cost, as depicted in Figure~\ref{figconfigs}. In Layout 1, we optimize inventory pooling by assigning all available vehicles to a central facility positioned at the midpoint of the grid. Conversely, in Layout 10, with 49 facilities, each facility is paired with one vehicle distributed across all blocks of the location grid. In this scenario, there is no benefit from inventory pooling. The other scenarios, which fall between these two edge cases, reflect varying effects of inventory pooling based on the number of open facilities and their distribution across the grid. Figure~\ref{figconfigstotalcost} illustrates the normalized total cost of the system across various experiments. Table~\ref{table:0} lists the parameters employed for these scenarios. In each experiment, the customer patience radius and normalized total facility cost are varied relative to each other.

 \subsection{Centralization vs. Distribution: Costs and Patience Trade-offs}
\begin{figure}[ht!]
\centering
\begin{tabular}{ccc}
\includegraphics[width=0.3\textwidth]{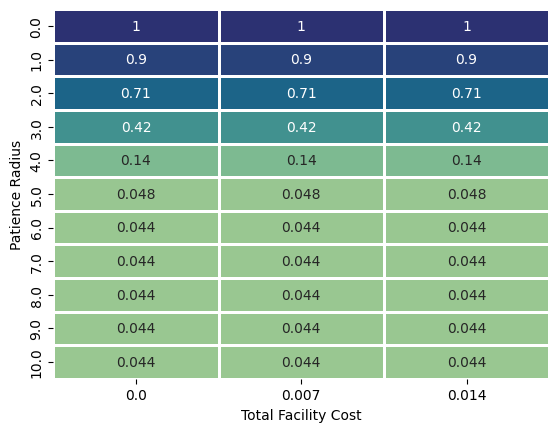}
& \includegraphics[width=0.3\textwidth]{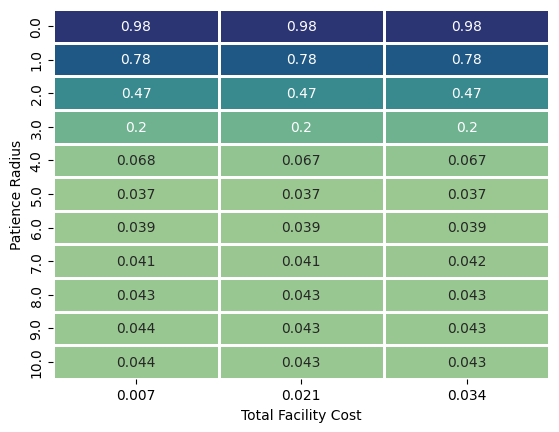}
&\includegraphics[width=0.3\textwidth]{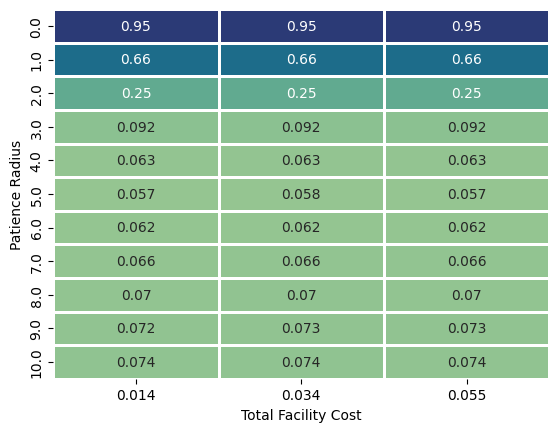} \\
Layout 1 &
Layout 2 &
Layout 3 \\
\includegraphics[width=0.3\textwidth]{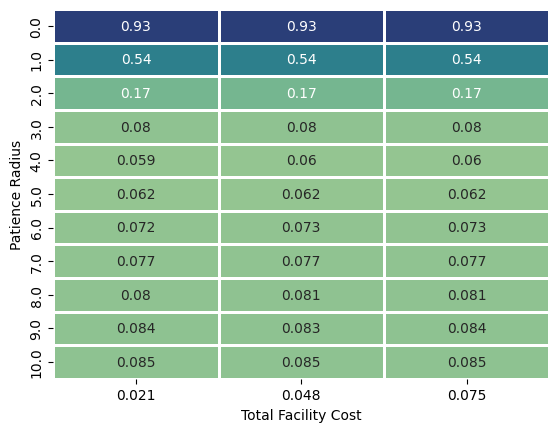} &
\includegraphics[width=0.3\textwidth]{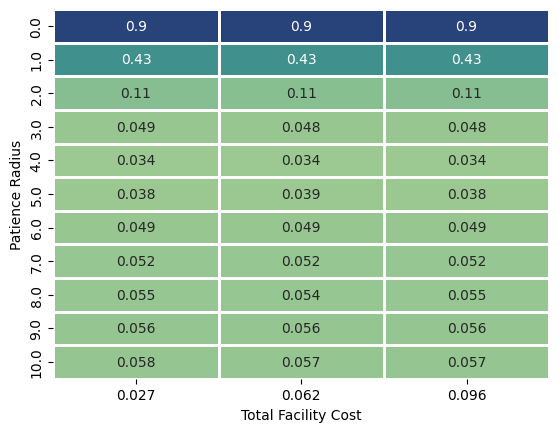} &
\includegraphics[width=0.3\textwidth]{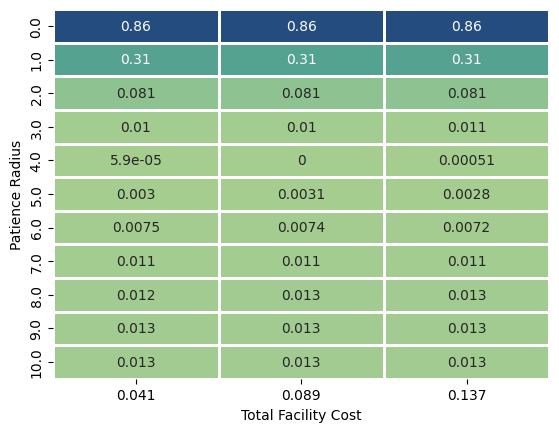}\\
Layout 4 &
Layout 5 &
Layout 6 \\
& & \\
\includegraphics[width=0.3\textwidth]{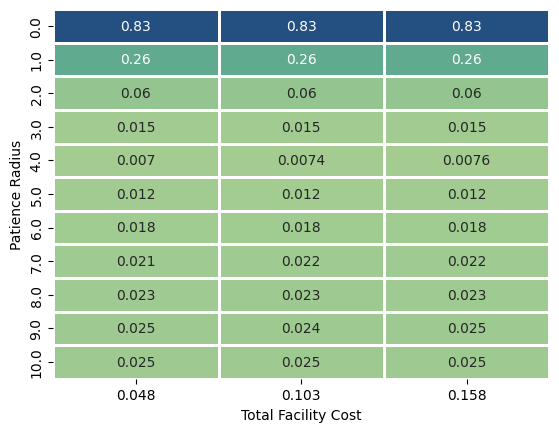} &
\includegraphics[width=0.3\textwidth]{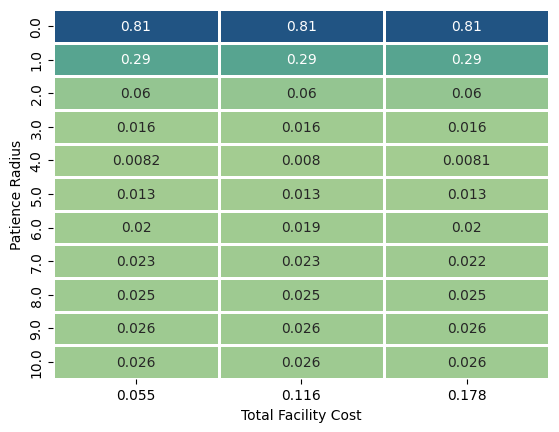} &
\includegraphics[width=0.3\textwidth]{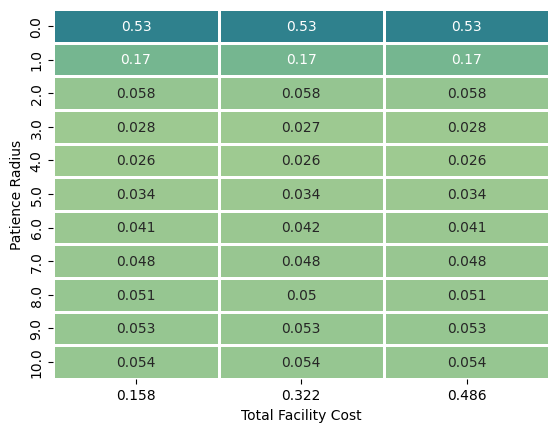} \\
Layout 7 &
Layout 8 &
Layout 9 \\
& & \\
\includegraphics[width=0.3\textwidth]{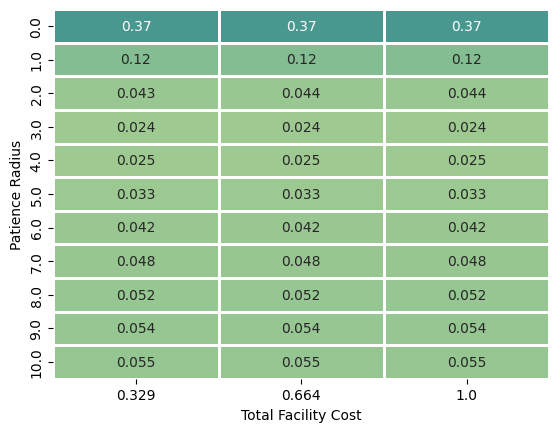} & & \\
Layout 10 & & \\
\end{tabular}
\caption{Normalized total lost sale cost of the system for different facility layouts}
\label{figconfigstotalcost1}
\end{figure}

 The results for Layout $1$, where inventory pooling is optimized, show that the total system cost decreases in this centralized layout.
However, moving to layouts with more facilities, with Layout $10$ having the highest number of open facilities, we observe that the total system cost rises drastically with increased facilities and patience levels. This can be attributed to the fact that when customers exhibit greater patience, opening numerous facilities becomes sub-optimal as it escalates the fixed cost of the system. Maintaining many open facilities is not a viable long-term strategy in systems with patient customers. Capitalizing on fewer open facilities is advantageous by accepting increased travel time from customers. Also, Figure~\ref{figconfigstotalcost1} depicts the system's normalized total lost sale cost in different experiments. Each experiment has variations in the customer patience radius and normalized total facility cost relative to each other. The findings indicate that maintaining more open facilities in systems with fewer patient customers is advisable for minimizing the total lost sale cost. This is more visible in Layout $10$ as shown in Figure~\ref{figconfigstotalcost1}, where the total lost sale cost reaches its minimum.

\begin{table}[ht!]
\centering
\caption{Parameters for random test scenarios}\label{table:0}
%\begin{adjustbox}{width=0.6\columnwidth,center}
\begin{tabular}{cc }
 \hline
 Parameter & Value\\ 
 \hline
 Dimension of location grid & $7\times 7$\\ 
 Total number of available vehicles  & 49\\
 Cost of keeping a facility open per unit of time& \$5, \$10,\$15 \\
 Distance travel cost per block traveled& \$5 \\
 Lost sale cost per unit of lost demand & \$5 \\
 Customer patience radius & 1 to 10 blocks away\\
 Total number of customers & 100\\
 Time unit horizon & A day \\
 \hline
\end{tabular}
%\end{adjustbox}
\end{table}

To assess a scenario closer to real-world conditions, we construct a situation with no associated cost for keeping a facility open during each time unit. This assumption aligns with the possibility of AVs being stationed at street parking without incurring any cost in practical situations. The patience radius is neither very low nor high in typical situations involving average customers. We opt for a moderate value of $3$ blocks as the patience radius. Another realistic assumption is that the lost sale cost per unit is significantly higher in the system than the travel cost per block. Table~\ref{table:1} illustrates the parameters employed for this scenario.

\begin{table}[ht!]
\centering
\caption{Parameters for real case test scenarios}\label{table:1}
%\begin{adjustbox}{width=0.6\textwidth,center}
\begin{tabular}{ ccc }
 \hline
 Parameter & Value\\ 
 \hline
 Dimension of location grid & 7X7\\ 
 Total number of available vehicles  & 49\\
 Cost of keeping a facility open per unit of time& \$0\\
 Distance travel cost per block traveled& \$1 \\
 Lost sale cost per unit of lost demand & \$50 \\
 Customer patience radius & Up to 3 blocks away\\
 Total number of customers & 100\\
 Time unit horizon & A day\\
 \hline
\end{tabular}
%\end{adjustbox}
\end{table}

\begin{figure}[ht!]
   \centering
       \includegraphics[scale=0.5]{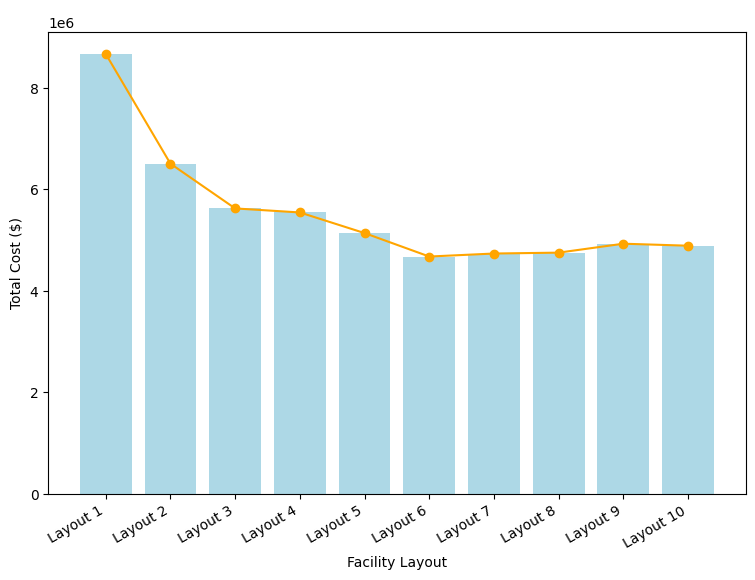}
       \caption{Average total cost of the system for different facility configurations for 100 independent days}
       \label{fig:1}
   \end{figure}

As depicted in Figure \ref{fig:1}, the total system cost demonstrates a decrease as the number of facilities increases. Notably, there is evidence of diminishing returns when expanding the number of facilities. Beyond nine facilities, the overall cost remains relatively stable. This phenomenon can be attributed to the absence of fixed costs for keeping facilities open. Given that the customers are regular and there is no urgency for additional facilities, the marginal benefit of adding more facilities diminishes. Instead, it increases the total system cost without significantly improving demand satisfaction.

\begin{figure}[ht!]

   \centering
       \includegraphics[scale=0.5]{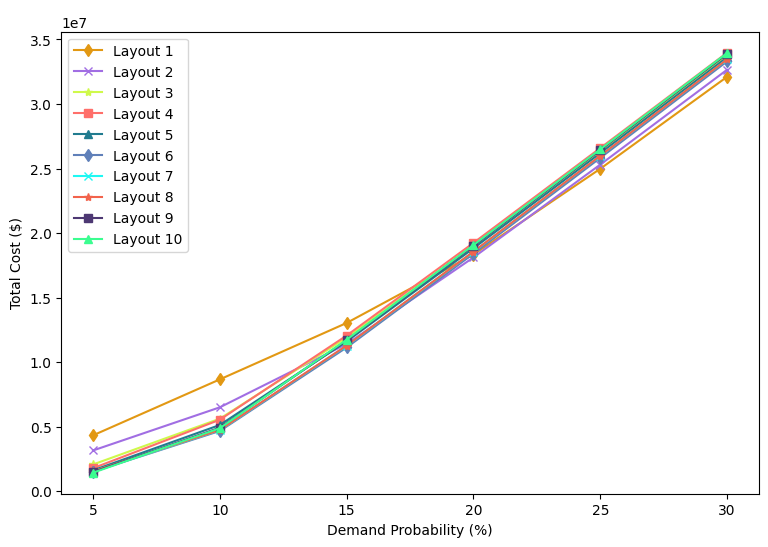}
       \caption{Average total cost of the system with varying demand rates  from $5$\% to $30$\% for different facility layouts for 100 independent days}
       \label{fig:2}
   \end{figure}

\subsection{Sensitivity Analysis of AV Routing: Demand and Patience Impact}   
\begin{figure}[ht!]
   \centering
       \includegraphics[scale=0.5]{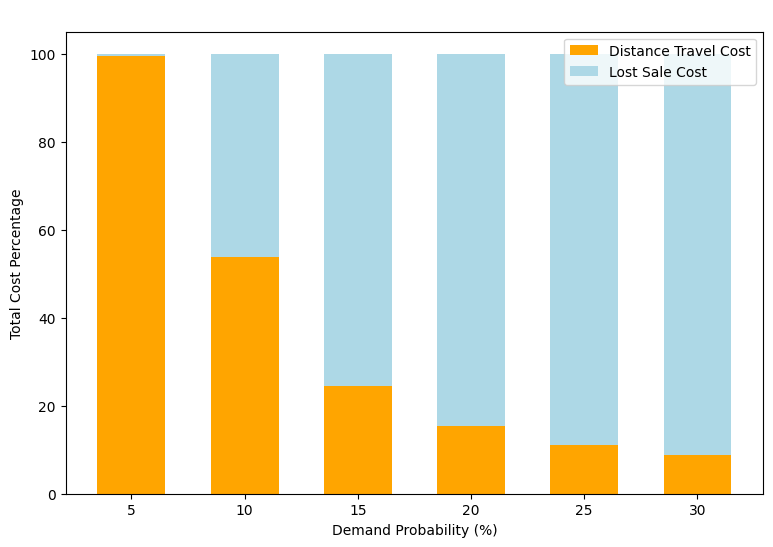}
       \caption{Average total travel and lost sale cost of the system with varying demand rates  from $5$\% to $30$\% for different facility layouts for 100 independent days}
       \label{fig:3}
   \end{figure}

A sensitivity analysis is conducted on the demand rate, as illustrated in Figure \ref{fig:2}. The probability of customer presence at each block is assumed to be a Bernoulli distribution, with the success probability varied from $5$\% to $30$\%. Other parameters remained constant at their previously specified values in Table~\ref{table:1}. As shown in Figure~\ref{fig:2}, the system's total cost increases with increased success probability rate, which is anticipated given that a higher number of customers results in increased travel costs and lost sales costs. Interestingly, the scenario with a $1$ facility demonstrates linear growth, whereas the others exhibit exponential growth. This shift is related to the percentage of lost sales costs. In other words, as demand rises, a larger proportion of the total cost comprises lost sales costs. Consequently, the difference between the number of facilities becomes less visible as demand increases, as all vehicles are immediately occupied, and the lost sales costs remain roughly consistent. Figure \ref{fig:3} visually represents this concept, illustrating how lost sales costs become a significantly more significant fraction of the total cost as the demand rate rises.

\begin{figure}[ht!]
   \centering
       \includegraphics[scale=0.5]{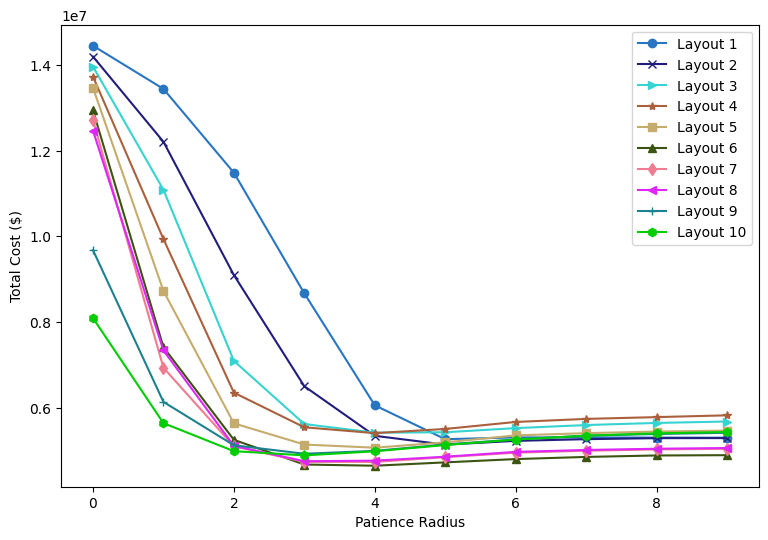}
       \caption{Average total cost of the system with varying customer patience level from 0 to 10 for different facility layouts and 100 independent days}
       \label{fig:4}
   \end{figure}

A crucial parameter chosen for sensitivity analysis is customers' patience level. In the initial trial, this limit was set up to 3 blocks from the customer. To evaluate the effect of the patience level of customers on the total cost of the system, we adjust this parameter across a range from $0$ to $10$ for all customers and present the outcomes in Figure \ref{fig:4}.
The results show that as customers exhibit increased patience, the system's total cost decreases for all facility layouts. However, an interesting finding is that the utterly centralized case, where we have only one open facility, initially proves to be the most expensive for very impatient customers but eventually emerges as the most cost-effective alternative as customer patience increases. This is a primary result of this paper, highlighting that the advantages of inventory pooling for stationing AVs amplify with customers' rising patience.

\begin{figure}[ht!]
   \centering
       \includegraphics[scale=0.5]{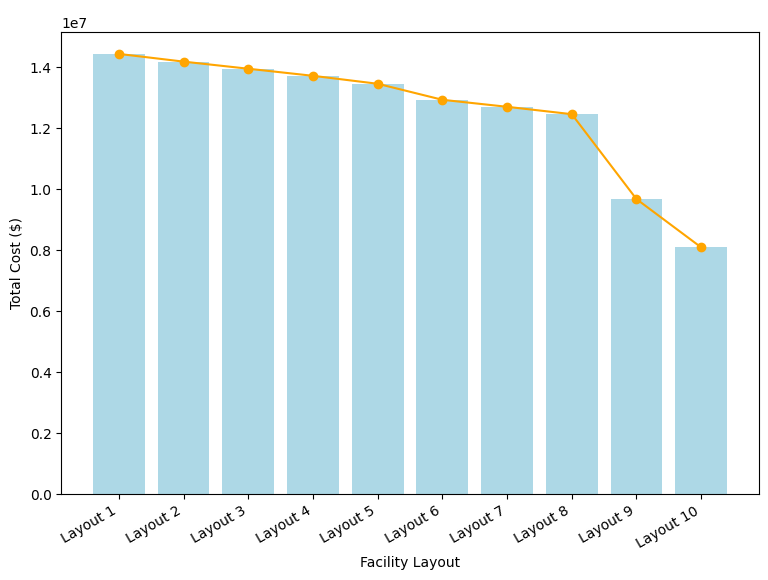}
       \caption{Average total cost of the system with customer patience level $= 0$,  for different facility configurations and 100 independent days}
       \label{fig:5}
   \end{figure}
   
In other words, maximizing the number of operational facilities is advantageous for customers with extremely low patience levels, each with a patience radius zero. 
Intuitively, since customers only accept a vehicle if it is currently available in their location (given that the patient radius is 0), it is most beneficial for the company to distribute AVs as widely as possible. Figure \ref{fig:5} visually confirms this observation, as the system's total cost consistently decreases with more facilities.

\begin{figure}[ht!]
   \centering
       \includegraphics[scale=0.5]{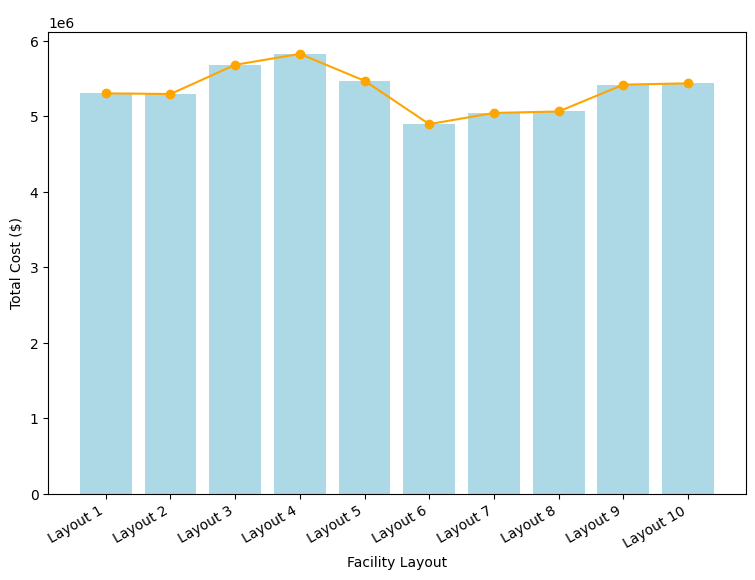}
       \caption{Average total cost of the system with customer patience level $= 10$,  for different facility configurations and 100 independent days}
       \label{fig:6}
   \end{figure}

Figure \ref{fig:6} presents a graph similar to Figure \ref{fig:5}, but with customers displaying more remarkable patience, having a waiting radius of 10 blocks. Unlike the linear relationship observed with a patience radius of 0, the graph exhibits a spike between 2 and 4 facilities, a decline for 5 through 7 facilities, and a rise for more facilities. This spike is attributed to the assignment of vehicles to customers. The proposed heuristic relies solely on the current availability of vehicles to match them with customers at each time unit. For scenarios with four facilities, the customer's nearest vehicle might be far away if the other facilities have no available vehicles. Given the high patience of customers in this case, they are willing to wait for the distant vehicle to travel, whether it is 8, 9, or 10 blocks. However, this incurs substantial travel costs for the vehicle, which covers a long distance. Moreover, it is conceivable that while the distant customer is being served, another customer could arrive closer to the vehicle's initial facility. In such a situation, a lost sale cost may be incurred if the vehicle undertaking the long-distance trip is the last available one and no more customers are present. In this context, it would have been more strategic for the company to reject the initial customer, bear the lost sales cost, and reserve that vehicle for a closer customer in the future. Essentially, the platform should be able to anticipate and decide when to strategically decline a customer to preserve a vehicle for subsequent customers. This decision could depend on the distance to the customer and the number of vehicles returning. Unlike other business models, each vehicle's travel duration and return times are known. Hence, assigning a vehicle to a more distant customer might be worthwhile if numerous vehicles are anticipated to return to a facility soon. On the other hand, if few vehicles are in the facility and few are expected to return soon, saving the vehicles for the future, closer customers would be more cost-effective.

\section{Discussion and conclusion} \label{conclusion}

This paper introduced a novel mathematical model for allocating AVs to customers, taking into account their patience levels to minimize the overall system cost. A heuristic algorithm was employed to solve the proposed model and enhance practicality.
Several insightful conclusions emerged from this study. Firstly, the effectiveness of adding more facilities diminished as the existing number of facilities increased. Additionally, the total cost of the system increased with the demand rate. As demand rose, the total cost for varying facility numbers began to converge, primarily due to the growing proportion of costs attributed to lost sales, which became similar across all facility configurations. 
Moreover, the investigation demonstrated that with increased customer patience, the advantages of inventory pooling escalated, particularly in scenarios where no fixed cost is associated with keeping the facility open during the time interval. While an optimal policy was not identified due to the model's complexity, the proposed heuristic yielded several noteworthy insights. This paper also highlighted the need for an optimal policy to consider current customers and available vehicles, future returning vehicles, and potential future demand. 

Furthermore, this paper opens avenues for future research in different ways. First, a personalized patience radius could be explored, acknowledging that individuals possess different patience levels that may vary under various circumstances. Separating customers into classes with distinct patience levels could better reflect the diverse market. Another potential area for exploration involves introducing capacity constraints on vehicles, allowing them to pick up multiple customers. The rigid policy of vehicles returning to their facility before becoming available could be relaxed, with alternative vehicles idling after dropping off a customer or returning to the nearest facility with available space. These policy adjustments may reduce travel costs and merit further investigation.

Additionally, changes to the demand pattern could be considered. All locations are treated equally, but incorporating regional popularity differences or time-dependent demand variations (e.g., morning suburban to downtown and evening downtown to suburban commutes) could introduce seasonality and impact facility and vehicle placement decisions. While this study did not yield an exact solution, it uncovered fundamental patterns that could enhance heuristic policies for decision-makers in geospatially-driven demand scenarios. The research also spawned intriguing avenues for further exploration, posing additional research questions and potential applications. Extending the problem to a two-echelon delivery system \cite{moradi2024two} and considering same-day delivery concepts \cite{moradi2025prize} and \cite{boroujeni2025last} could be promising research topics.

\section*{CRediT authorship contribution statement}
\textbf{Niloufar Mirzavand Boroujeni:} Conceptualization, Data curation, Formal analysis, Investigation, Methodology, Software, Supervision, Validation, Visualization, Roles/Writing - original draft, and Writing - review \& editing.
\textbf{Nasim Mirzavand Boroujeni:} Formal analysis, Resources, Software, Validation, Roles/Writing - original draft; and Writing - review \& editing.
\textbf{Nima Moradi \&  Saeed Jamalzadeh:} Data curation, Formal analysis, Resources, Software, Validation, Visualization, and Writing - review \& editing.

\section*{Conflict of interest} 
The authors have no conflicts of interest to disclose.

\section*{Funding} The authors received no financial support for this paper's research, authorship, and publication.

\section*{Data Availability Statement}
Data, models, and codes are available upon request.

%\clearpage

\bibliography{sn-bibliography}% common bib file
%% if required, the content of .bbl file can be included here once bbl is generated
%%\input sn-article.bbl

\end{document}